\documentclass[reqno,a4paper,final,12pt]{amsart}
\usepackage[latin1]{inputenc}
\usepackage{amssymb,amstext,amsthm,eucal}
\usepackage{color}
\usepackage{array}
\usepackage{a4wide}
\usepackage[notref,notcite]{showkeys}
\usepackage[all]{xy}
\usepackage{hyperref}
\hypersetup{%
  pdftitle   = {A geometric construction of the exceptional Lie
    algebras F4 and E8},
  pdfkeywords = {supergravity, exceptional Lie algebra, E8, Killing
    superalgebra, Killing spinors},
  pdfauthor  = {José Figueroa-O'Farrill},
  pdfcreator = {\LaTeX\ with package \flqq hyperref\frqq}
}
\newcommand{\half}{\tfrac12}
\newcommand{\be}{\boldsymbol{e}}

\newcommand{\bomega}{\boldsymbol{\omega}}
\newcommand{\bv}{\boldsymbol{v}}
\newcommand{\fe}{\mathfrak{e}}
\newcommand{\ff}{\mathfrak{f}}

\newcommand{\fk}{\mathfrak{k}}

\newcommand{\fM}{\mathfrak{M}}
\newcommand{\fS}{\mathfrak{S}}
\newcommand{\fso}{\mathfrak{so}}

\newcommand{\Cl}{\mathrm{C}\ell}
\renewcommand{\1}{\mathbf{1}}

\newcommand{\Spin}{\mathrm{Spin}}
\newcommand{\SO}{\mathrm{SO}}
\renewcommand{\O}{\mathrm{O}}
\ifx\Sp\UnDeFiNeD
\newcommand{\Sp}{\mathrm{Sp}}
\else
\renewcommand{\Sp}{\mathrm{Sp}}
\fi

\newcommand{\RR}{\mathbb{R}}
\newcommand{\CC}{\mathbb{C}}
\newcommand{\KK}{\mathbb{K}}
\newcommand{\HH}{\mathbb{H}}

\newcommand{\ZZ}{\mathbb{Z}}
\newcommand{\eL}{\mathcal{L}}

\DeclareMathOperator{\End}{End}
\DeclareMathOperator{\Hom}{Hom}

\DeclareMathOperator{\id}{id}

\newcommand{\MUNCH}[1]{\relax}

\allowdisplaybreaks[1]
\begin{document}
\title[Geometric construction of exceptional Lie algebras]{A geometric
  construction of the exceptional Lie algebras $F_4$ and $E_8$}
\author{José Figueroa-O'Farrill}
\address{Maxwell Institute and School of Mathematics, University of
  Edinburgh, UK}
\email{j.m.figueroa@ed.ac.uk}
\thanks{EMPG-07-12}
\begin{abstract}
  We present a geometric construction of the exceptional Lie algebras
  $F_4$ and $E_8$ starting from the round spheres $S^8$ and $S^{15}$,
  respectively, inspired by the construction of the Killing
  superalgebra of a supersymmetric supergravity background.
\end{abstract}
\maketitle

\section{Introduction}
\label{sec:introduction}

The Killing--Cartan classification of simple Lie algebras over the
complex numbers is well known: there are four infinite families
$A_{n\geq 1}$, $B_{n\geq 2}$, $C_{n \geq 3}$ and $D_{n \geq 4}$, with
the range of ranks chosen to avoid any overlaps, and five exceptional
cases $G_2$, $F_4$, $E_6$, $E_7$ and $E_8$.  Whereas the classical
series ($A$-$D$) correspond to matrix Lie algebras, and indeed their
compact real forms are the Lie algebras of the special unitary groups
over $\CC$ ($A$), $\HH$ ($B$) and $\RR$ ($C$ and $D$), the exceptional
series do not have such classical descriptions; although they can be
understood in terms of more exotic algebraic structures such as
octonions and Jordan algebras.  There is, however, a uniform
construction of all exceptional Lie algebras (except for $G_2$) using
spin groups and their spinor representations, described in Adams'
posthumous notes on exceptional Lie groups \cite{AdamsExceptional}
and, for the special case of $E_8$, also in \cite{GSW}.  This
construction, once suitably geometrised, is very familiar to
practitioners of supergravity.  The purpose of this note is to present
this geometrisation, perhaps as an invitation for differential
geometers to think about supergravity.

Indeed in supergravity there is a geometric construction which
associates a Lie superalgebra to any supersymmetric supergravity
background: typically a lorentzian spin manifold with extra geometric
data and with a notion of privileged spinor fields, called Killing
spinors.  The resulting superalgebra is called the Killing
superalgebra because it is constructed out of these Killing spinors
and Killing vectors.  The Killing superalgebra for general ten- and
eleven-dimensional supergravities is constructed in \cite{FMPHom,
  EHJGMHom}.  In this note we will apply this construction not to
supergravity backgrounds, but to riemannian manifolds without any
additional structure.  The relevant notion of Killing spinor is then
that of a \emph{geometric Killing spinor}: a nonzero section
$\varepsilon$ of the spinor bundle satisfying
\begin{equation*}
  \nabla_X \varepsilon = \half X \cdot \varepsilon~,
\end{equation*}
where $X$ is any vector field and the dot means the Clifford action.
We will apply this construction to the unit spheres $S^7 \subset
\RR^8$, $S^8 \subset \RR^9$ and $S^{15} \subset \RR^{16}$ and in this
way obtain the compact real Lie algebras $\fso_9$, $\ff_4$ and
$\fe_8$, respectively.  It is curious that these three spheres are
linked by the exceptional Hopf fibration which defines the octonionic
projective line,
\begin{equation*}
  \xymatrix{ S^7 \ar[r] & S^{15} \ar[d] \\
    & S^8 }
\end{equation*}
and it is natural to wonder whether their Killing superalgebras are
similarly related.  We will not answer this question here.

This note is organised as follows.  In Section~\ref{sec:algebra} we
briefly review the relevant notions of Clifford algebras, spin groups
and their spinorial representations.  In Section~\ref{sec:KSA} we
define the Killing superalgebra after introducing the basic notions of
Killing spinors and Bär's cone construction.  In
Section~\ref{sec:S7S8S15} we construct the Killing superalgebras of
the round spheres $S^7$, $S^8$ and $S^{15}$ and show that they are
isomorphic to the compact real Lie algebras $\fso_9$, $\ff_4$ and
$\fe_8$, respectively.  Finally in Section~\ref{sec:conc} we discuss
some open questions motivated by the results presented here.

\section{Spinorial algebra}
\label{sec:algebra}

In this section we start with some algebraic preliminaries on
euclidean Clifford algebras and spinors in order to set the notation.
We will be sketchy, but fuller treatments can be found, for example,
in \cite{ABS,LM,Harvey}.

\subsection{Clifford algebras and Clifford modules}
\label{sec:clifford}

Let $V$ be a finite-dimensional real vector space with a
positive-definite euclidean inner product $\left<-,-\right>$.
The \textbf{Clifford algebra} $\Cl(V)$ is the associated algebra with
unit generated by $V$ and the identity $\1$ subject to the Clifford
relations
\begin{equation}
  \label{eq:clalg}
  \bv^2 = - \left<\bv,\bv\right> \1
\end{equation}
for all $\bv \in V$.  More formally, the Clifford algebra is the
quotient of the tensor algebra of $V$ by the two-sided ideal generated
by the Clifford relations.  Since the Clifford relations---having terms
of degree 0 and degree 2---are not homogeneous in the natural grading
of the tensor algebra, $\Cl(V)$ is not graded but only filtered.  The
associated graded algebra is the exterior algebra $\Lambda V$, to
which it is isomorphic as a vector space.  Nevertheless since the
terms in \eqref{eq:clalg} have even degree, $\Cl(V)$ is
$\ZZ_2$ graded
\begin{equation}
  \Cl(V) = \Cl(V)_0 \oplus \Cl(V)_1~,
\end{equation}
with vector-space isomorphisms $\Cl(V)_0 \cong \Lambda^{\text{even}}V$
and $\Cl(V)_1 \cong \Lambda^{\text{odd}}V$.  These isomorphisms can be
seen explicitly as follows.  Relative to an orthonormal basis $\be_i$ for
$V$, the Clifford relations become
\begin{equation}
  \label{eq:clifford2}
  \be_i \be_j + \be_j \be_i = - 2 \delta_{ij} \1~,
\end{equation}
which shows that up to terms in lower order we may always
antisymmetrise any product $\be_{i_1} \be_{i_2} \dots \be_{i_k}$ in
$\Cl(V)$ without ever changing the parity.  The Clifford algebra of
$\RR^n$ generated by $\1$ and $\be_i$ subject to \eqref{eq:clifford2}
is denoted $\Cl_n$.  As a real associative algebra with unit it is
isomorphic to one or two copies of matrix algebras, as shown in
Table~\ref{tab:Cln} for $n\leq 7$.  The higher values of $n$ are
obtained by Bott periodicity $\Cl_{n+8} \cong \Cl_n \otimes \RR(16)$,
where $\RR(16)$ is the algebra of $16 \times 16$ real matrices.

\begin{table}[h!]
  \centering
  \renewcommand{\arraystretch}{1.3}
  \begin{tabular}{|>{$}c<{$}|*{8}{>{$}c<{$}|}}
    \hline
    n & 0 & 1 & 2 & 3 & 4 & 5 & 6 & 7 \\
    \hline
    \Cl_n & \RR & \CC & \HH & \HH \oplus \HH & \HH(2) & \CC(4) &
    \RR(8) & \RR(8) \oplus \RR(8)\\
    [3pt]
    \hline
  \end{tabular}
  \vspace{8pt}
  \caption{Clifford algebras $\Cl_n$, where $\KK(m)$ denotes the
    algebra of $m\times m$ matrices with entries in $\KK$.}
  \label{tab:Cln}
\end{table}

Since matrix algebras have a unique irreducible representation (up to
isomorphism), we can easily read off the irreducible representations
of $\Cl_n$ from the table.  We see, in particular, that if $n$ is even
there is a unique irreducible representation, which is real for $n
\equiv 0,6 \pmod 8$ or quaternionic when $n \equiv 2,4 \pmod 8$;
whereas if $n$ is odd there are two inequivalent irreducible
representations, which are real when $n\equiv 7 \pmod 8$ and
quaternionic when $n\equiv 3 \pmod 8$, and form a complex conjugate
pair for $n \equiv 1 \pmod 4$.  These two inequivalent Clifford
modules are distinguished by the action of $\bomega :=
\be_1\be_2\cdots\be_n$, which for $n$ odd is central in $\Cl_n$.  This
element obeys $\bomega^2 = (-1)^{n(n+1)/2} \1$, whence it is a complex
structure for $n \equiv 1 \pmod 4$, in agreement with the table.  The
dimension of one such irreducible \textbf{Clifford module}, relative
to either $\RR$ if real or $\CC$ if not, is $2^{\lfloor n/2\rfloor}$.
We will use the notation $\fM$ for the unique irreducible Clifford
module in even dimension, $\fM_\pm$ for the irreducible Clifford
modules for $n\equiv 3 \pmod 4$.  For $n \equiv 1 \pmod 4$ we will let
$\fM$ denote the irreducible Clifford module on which $\bomega$ acts
like $+i$ and let $\overline\fM$ denote the irreducible module on
which $\bomega$ acts like $-i$.

\subsection{The spin group and spinor modules}
\label{sec:spin}

The Clifford algebra $\Cl(V)$ admits a natural Lie algebra structure
via the Clifford commutator.  The map $\Lambda^2V \to \Cl(V)$ given by
\begin{equation}
  \be_i \wedge \be_j \mapsto -\half \be_i \be_j~,
\end{equation}
for $i<j$, induces a Lie algebra homomorphism $\rho: \fso(V) \to
\Cl(V)$.  Moreover the action of $\fso(V)$ on $V$ is realised by the
Clifford commutator, so that if $A \in \fso(V)$ and $\bv \in V$, then
\begin{equation}
  A(\bv) = \rho(A) \bv - \bv \rho(A) \in \Cl(V)~.
\end{equation}
Exponentiating the image of $\rho$ in $\Cl(V)$ we obtain a connected
Lie group called $\Spin(V)$.  The subspace $V \subset \Cl(V)$ is
closed under conjugation by $\Spin(V)$ whence we obtain a map
$\Spin(V) \to \SO(V)$, whose kernel is the central subgroup consisting
of $\pm \1$.

Restricting an irreducible Clifford module $\fM$ (or $\fM_\pm$) to
$\Spin(V)$ we obtain a \textbf{spinor module}, which may or may not
remain irreducible.  Since $\Spin_n \subset \left(\Cl_n\right)_0 \cong
\Cl_{n-1}$, we can immediately infer the type of spinor module from
Table~\ref{tab:Cln}.  If $n\equiv 1\pmod 8$, $\fM \cong \fS \otimes
\CC$ is the complexification of the unique irreducible real spinor
module $\fS$, whereas if $n\equiv 5 \pmod 8$, $\fM \cong \fS$, but
$\fS$ possesses a $\Spin_n$-invariant quaternionic structure, whence
$\overline\fM \cong \fS$ as well.  For $n\equiv 3 \pmod 4$, $\fM_\pm
\cong \fS$.  For odd $n$, the spinor module is real if $n\equiv 1,7
\pmod 8$ and quaternionic otherwise and its dimension (over $\RR$ if
real and over $\CC$ otherwise) is again $2^{(n-1)/2}$.  For even $n$,
the unique irreducible Clifford module decomposes (perhaps after
complexification) into two inequivalent $\Spin_n$ modules, called
\textbf{half-} (or \textbf{chiral}) \textbf{spinor} modules.  They are
denoted $\fS_\pm$ if $n\equiv 0 \pmod 4$ and $\fS$ and $\overline\fS$
if $n \equiv 2 \pmod 4$.  They are real if $n\equiv 0 \pmod 8$,
quaternionic if $n \equiv 4 \pmod 8$ and complex otherwise.  If
$n\equiv 6 \pmod 8$ then it is the complexification of $\fM$ which
decomposes $\fM \otimes \CC \cong \fS \oplus \overline \fS$.  In all
cases, the dimension, computed relative to the appropriate field for
the type, is $2^{(n-2)/2}$.

\subsection{Spinor inner products}
\label{sec:inner}

The Clifford algebra $\Cl(V)$ has a natural antiautomorphism defined
by $-\id_V$ on $V$.  On a given irreducible Clifford module $\fM$ (or
$\fM_\pm$) there always exists an inner product $\left(-,-\right)$
which realises this automorphism; that is, such that
\begin{equation}
  \label{eq:ClIP}
  \left(\bv \cdot \varepsilon_1, \varepsilon_2\right) = 
  - \left(\varepsilon_1, \bv \cdot \varepsilon_2\right)~,
\end{equation}
for all $\bv \in V$ and $\varepsilon_i \in \fM$.  It follows that
$\left(-,-\right)$ is $\Spin(V)$-invariant; indeed,
\begin{equation}
  \left(\be_i\be_j \cdot \varepsilon_1, \varepsilon_2\right) = 
  - \left(\varepsilon_1, \be_i\be_j \cdot \varepsilon_2\right)~.
\end{equation}
In positive-definite signature, $\left(-,-\right)$ is either symmetric
or hermitian, depending on the type of representation, and
positive-definite \cite{Harvey}.

The Clifford action $V \otimes \fM \to \fM$ induces a map,
suggestively denoted $[-,-]: \fM \otimes \fM \to V$, via the above
inner product on $\fM$ and the euclidean inner product
$\left<-,-\right>$ on $V$.  Explicitly, we have that for all $\bv \in
V$ and $\varepsilon_i\in\fM$,
\begin{equation}
  \label{eq:bispinor1}
  \left<[\varepsilon_1,\varepsilon_2],\bv\right> = \left(
    \varepsilon_1, \bv \cdot \varepsilon_2\right)~.
\end{equation}

\section{The Killing superalgebra}
\label{sec:KSA}

In this section we will define the Killing superalgebra of a
riemannian spin manifold admitting Killing spinors.

\subsection{Spin manifolds}
\label{sec:spin-manifolds}

Let $(M,g)$ be an $n$-dimensional riemannian manifold and let $\O(M)$
denote the bundle of orthonormal frames. It is a principal
$\O_n$-bundle.  If the manifold is orientable, we can restrict
ourselves consistently to oriented orthonormal frames.  In this case,
the subbundle $\SO(M)$ of oriented orthonormal frames is a principal
$\SO_n$-bundle.  The obstruction to orientability is measured by the
first Stiefel--Whitney class $w_1 \in H^1(M;\ZZ_2)$.  If $(M,g)$ is
orientable one can ask whether there is a principal $\Spin_n$-bundle
$\Spin(M)$ lifting the oriented orthonormal frame bundle $\SO(M)$;
that is, admitting a bundle map $\Spin(M)\to\SO(M)$ covering the
identity and restricting fibrewise to the natural homomorphism
$\Spin_n \to \SO_n$.  The obstruction to the existence of such a lift
is measured by the second Stiefel--Whitney class $w_2 \in
H^2(M;\ZZ_2)$ and if it vanishes the manifold $(M,g)$ is said to be
spin.  Spin structures $\Spin(M)$ on $M$ need not be unique: they are
measured by $H^1(M;\ZZ_2) = \Hom(\pi_1M,\ZZ_2)$, which we can
understand as assigning a sign (consistently) to every noncontractible
loop.  In this section we will assume our manifolds to be spin and
that a choice of spin structure has been made.  The main examples in
this note are spheres, which are spin---indeed, the total space of the
spin bundle of $S^n$ is the Lie group $\Spin_{n+1}$---and, since they
are simply-connected, have a unique spin structure.

If $\fM$ is a $\Cl_n$-module, then it is also a (perhaps reducible)
$\Spin_n$-module and we may form the \textbf{spinor bundle}
\begin{equation*}
  S(M) := \Spin(M) \times_{\Spin_n} \fM
\end{equation*}
over $M$ as an associated vector bundle to the spin bundle.
Furthermore we have a fibrewise action of the Clifford bundle
$\Cl(TM)$ on $S(M)$.  The spinor inner products globalise to give
an inner product on $S(M)$.

The Levi-Cività connection on the orthonormal frame bundle of $(M,g)$
induces a connection on $\Spin(M)$ and hence on any associated vector
bundle.  In particular we have a \textbf{spin connection} on $S(M)$
and $\$(M)$.  This is a map
\begin{equation*}
  \nabla : \Gamma(S(M)) \to \Omega^1(M;S(M))~,
\end{equation*}
and similarly for $\$(M)$, and it allows us to write down
interesting equations on spinors.  One such equation is the Killing
spinor equation, which is the subject of the next section.  A classic
treatise on this equation is \cite{BFGK}.

\subsection{Killing spinors}
\label{sec:killing-spinors}

Throughout this section we will let $(M^n,g)$ be a spin manifold with
chosen spinor bundle $S(M)$ on which we have a fibrewise action of the
Clifford bundle $\Cl(TM)$ and a $\Spin_n$-invariant inner product
which in addition satisfies equation \eqref{eq:ClIP}.  A nonzero
$\varepsilon \in \Gamma(S(M))$ is said to be a \textbf{(real) Killing
  spinor} if for all vector fields $X$,
\begin{equation}
  \label{eq:KS}
  \nabla_X \varepsilon = \lambda X \cdot \varepsilon~,
\end{equation}
where $\lambda \in \RR$ is the \textbf{Killing constant}.  The origin
of the name is that if $\varepsilon_i$, $i=1,2$, are Killing spinors,
then the vector field $V:=[\varepsilon_1,\varepsilon_2]$ defined
by equation \eqref{eq:bispinor1} is a Killing vector.  Indeed, for all
vector fields $X,Y$,
\begin{align*}
  g(\nabla_XV,Y) &= \left(\nabla_X\varepsilon_1, Y \cdot
    \varepsilon_2\right) +   \left(\varepsilon_1, Y \cdot
    \nabla_X\varepsilon_2\right) & \tag{by definition of $\nabla$}\\
  &= \lambda \left(X \cdot \varepsilon_1, Y \cdot \varepsilon_2\right)
  + \lambda \left(\varepsilon_1, Y \cdot X \cdot
    \varepsilon_2\right) & \tag{using equation~\eqref{eq:KS}}\\
  &= - \lambda \left(\varepsilon_1, X \cdot Y \cdot
    \varepsilon_2\right) + \lambda \left(\varepsilon_1, Y \cdot X
    \cdot \varepsilon_2\right)~, & \tag{using equation
    \eqref{eq:ClIP}}
\end{align*}
which is manifestly skewsymmetric in $X,Y$, whence we conclude that
\begin{equation*}
  g(\nabla_XV,Y) + g(\nabla_Y V,X) = 0~,
\end{equation*}
which is one form of Killing's equation.

\subsection{The cone construction}
\label{sec:cone}

The problem of determining which riemannian manifolds admit real
Killing spinors was the subject of much research until it was
elegantly solved by Bär \cite{Baer} via the cone construction.  We
will assume that the Killing constant $\lambda$ has been set to $\pm
\half$ by rescaling the metric, if necessary.  Let $(\overline M,
\overline g)$ denote the \textbf{(deleted) cone} over $M$, defined by
$\overline M = \RR^+ \times M$ and $\overline g = dr^2 + r^2 g$, where
$r>0$ is the coordinate on $\RR^+$.  Bär observed that there is a
one-to-one correspondence between Killing spinors on $M$ and parallel
spinors on the cone $\overline M$.  More precisely, if $n = \dim M$ is
even, there is an isomorphism between Killing spinors on $M$ with
Killing constant $\pm \half$ and parallel spinors on $\overline M$;
the choice of sign having to do with the choice of embedding $\Cl_n
\subset \Cl_{n+1}$.  If on the other hand $n$ is odd, then the space
of Killing spinors on $M$ with Killing constant $\pm \half$ is
isomorphic to the space of parallel half-spinors on $\overline M$, the
chirality depending on the sign of the Killing constant.  Together
with a theorem of Gallot \cite{Gallot} which says that the cone of a
complete manifold is either flat or irreducible, the above observation
reduces the problem of determining the complete riemannian manifolds
admitting real Killing spinors to a holonomy problem which was solved
by Wang in \cite{Wang}.  If $(M,g)$ is not complete, its cone may be
reducible, but if so it can be shown to be locally a product of
subcones and applying Bär's results to each of the subcones allows one
to write local forms for the metrics on $M$ in terms of (double)
warped products \cite{FigLeiSim}.

For example, in the case of $M=S^n$, the cone is $\overline M =
\RR^{n+1}\setminus\{0\}$, but the metric extends smoothly to the
origin.  The space of parallel (half-)spinors on $\RR^{n+1}$ is
isomorphic to the relevant (half-)spinor representation of
$\Spin_{n+1}$.

\subsection{The Killing superalgebra}
\label{sec:killing-algebra}

To a riemannian manifold admitting real Killing spinors we may
associate an algebraic structure called the Killing superalgebra which
extends the Lie algebra of isometries in the following way.  The
underlying vector space is $\fk = \fk_0 \oplus \fk_1$ where $\fk_0$ is
the Lie algebra of isometries and $\fk_1$ is the space of Killing
spinors with $\lambda = \half$.  (There is a similar story for
$\lambda =-\half$.)  The bracket on $\fk$ consists of three pieces:
the Lie bracket on $\fk_0$, a map $\fk_0 \otimes \fk_1 \to \fk_1$ and
a map $\fk_1\otimes \fk_1 \to \fk_0$.  Depending on dimension and
signature, the latter map may be symmetric or antisymmetric, whence
the resulting bracket might correspond (if the Jacobi identity is
satisfied) to a Lie algebra or a Lie superalgebra.  In the riemannian
examples in this section we will recover Lie algebras, but in the
lorentzian examples common in supergravity the similar construction
leads to Lie superalgebras.  Let us now define these maps.

The map $\fk_1 \otimes \fk_1 \to \fk_0$ is induced from the algebraic
map $[-,-]$ in equation \eqref{eq:bispinor1}, which explains the
notation.  As we saw before the image indeed consist of Killing vector
fields.

The map $\fk_0 \otimes \fk_1 \to \fk_1$ is given by the
spinorial Lie derivative of Lichnerowicz and Kosmann(-Schwarzbach)
\cite{Kosmann} and which we now define.  If $X$ is a vector field on
$M$, then let $A_X : TM \to TM$ denote the endomorphism of the
tangent bundle defined by $A_XY = -\nabla_Y X$, for $\nabla$ the
Levi-Cività connection.  The vector field $X$ is Killing if and only
if $A_X$ is skewsymmetric relative to the metric; that is, if and only
if $A_X \in \fso(TM)$.  Let $\rho: \fso(TM) \to \End(S(M))$ denote the
spin representation and define the \textbf{spinorial Lie derivative}
along a Killing vector $X$ by
\begin{equation}
  \label{eq:spinLie}
  \eL_X = \nabla_X + \rho(A_X)~.
\end{equation}
In fact, this Lie derivative makes sense on sections of any vector
bundle associated to the orthonormal frame bundle provided that we
substitute $\rho$ for the relevant representation.  For instance, on
the tangent bundle itself, we have
\begin{equation*}
  \eL_X Y = \nabla_X Y + A_X Y = \nabla_X Y - \nabla_Y X = [X,Y]~,
\end{equation*}
as expected.  The spinorial Lie derivative satisfies the following
properties for all Killing vectors $X,Y$, spinors $\varepsilon$,
functions $f$ and arbitrary vector fields $Z$:
\begin{itemize}
\item $\eL_X$ is a derivation, so that
  \begin{equation}
    \label{eq:deriv}
    \eL_X(f \varepsilon) = X(f) \varepsilon + f \eL_X \varepsilon~;
  \end{equation}
\item $X \mapsto \eL_X$ is a representation of the Lie algebra of
  Killing vector fields:
  \begin{equation}
    \label{eq:morphism}
    \eL_X \eL_Y - \eL_Y \eL_X = \eL_{[X,Y]}~;
  \end{equation}
\item $\eL_X$ is compatible with Clifford multiplication:
  \begin{equation}
    \label{eq:clifford}
    \eL_X (Z \cdot \varepsilon) = [X,Z] \cdot \varepsilon + Z \cdot
    \eL_X \varepsilon~;
  \end{equation}
\item and $\eL_X$ preserves the Levi-Cività connection:
  \begin{equation}
    \label{eq:affine}
    \eL_X \nabla_Z - \nabla_Y \eL_X = \nabla_{[X,Z]}~.
  \end{equation}
\end{itemize}
It follows from equations \eqref{eq:clifford} and \eqref{eq:affine}
that the Lie derivative of a Killing spinor along a Killing vector is
again a Killing spinor.  Indeed, let $\varepsilon$ be a Killing spinor
and let $X$ be a Killing vector.  We have for all vector fields $Y$ that
\begin{align*}
  \nabla_Y \eL_X \varepsilon &= \eL_X \nabla_Y \varepsilon -
  \nabla_{[X,Y]} \varepsilon & \tag{using \eqref{eq:affine}}\\
  &= \lambda \eL_X \left(Y \cdot \varepsilon\right)- \lambda [X,Y]
  \cdot \varepsilon & \tag{since $\varepsilon$ is Killing}\\
  &= \lambda Y \cdot \eL_X \varepsilon~, & \tag{using
    \eqref{eq:clifford}}
\end{align*}
as advertised.  We define $[-,-]: \fk_0 \otimes \fk_1 \to \fk_1$ by
$[X,\varepsilon] := \eL_X \varepsilon$.

Of course, the existence of a bracket is not enough to conclude that
$\fk$ is Lie (super)algebra: one must also check the Jacobi identity.
The Jacobi identity is the vanishing of a tensor in $\fk\otimes
\Lambda^3\fk^*$.  Since $\fk = \fk_0 \oplus \fk_1$ and the bracket
respects the $\ZZ_2$ grading, there are four components to the Jacobi
identity.  The component in $\fk_0 \otimes \Lambda^3\fk_0$
vanishes due to the Jacobi identity of the Lie algebra $\fk_0$.  The
component in $\fk_1 \otimes \Lambda^2\fk_0^* \otimes\fk_1^*$ vanishes
because of the fact that $\fk_1$ is a representation of $\fk_0$;
indeed, this identity says that if $X,Y \in \fk_0$ and $\varepsilon
\in \fk_1$, then
\begin{equation*}
  [X,[Y,\varepsilon]] - [Y,[X,\varepsilon]] = [[X,Y],\varepsilon]~,
\end{equation*}
which is precisely equation (\ref{eq:morphism}).  The component in
$\fk_0 \otimes \Lambda^2\fk_1^* \otimes\fk_0^*$ vanishes because the
bracket $\fk_1 \otimes \fk_1 \to \fk_0$ is $\fk_0$-equivariant.
Indeed, if $X\in\fk_0$ and $\varepsilon_i\in\fk_1$ for $i=1,2$, then
for all vector fields $Y$,
\begin{align*}
  g\left([X,[\varepsilon_1,\varepsilon_2]],Y \right) &=
  g\left(\eL_X[\varepsilon_1,\varepsilon_2],Y \right)\\
  &= X g\left([\varepsilon_1,\varepsilon_2],Y \right) - 
  g\left([\varepsilon_1,\varepsilon_2],\eL_XY \right) & \tag{since $X$
    is Killing}\\
  &= X \left( \varepsilon_1, Y \cdot \varepsilon_2 \right) -
  \left( \varepsilon_1, \eL_XY \cdot \varepsilon_2 \right)\\
  &= \left( \eL_X \varepsilon_1, Y \cdot \varepsilon_2 \right) +
  \left( \varepsilon_1, \eL_X(Y \cdot \varepsilon_2) \right) -
  \left( \varepsilon_1, \eL_XY \cdot \varepsilon_2 \right)\\
  &= \left( \eL_X \varepsilon_1, Y \cdot \varepsilon_2 \right) +
  \left( \varepsilon_1, Y \cdot \eL_X\varepsilon_2\right) & \tag{using
  (\ref{eq:clifford})}\\
  &= g\left([[X,\varepsilon_1],\varepsilon_2],Y \right) +
  g\left([\varepsilon_1,[X,\varepsilon_2]],Y \right)~.
\end{align*}
The final component of the Jacobi identity lives in the
$\fk_0$-invariant subspace of $\fk_1 \otimes \Lambda^3\fk_1^*$.  This
identity does not seem to follow formally from the construction, but
requires a case-by-case argument.  In some cases it follows
because there simply are no $\fk_0$-invariant tensors in $\fk_1
\otimes \Lambda^3\fk_1^*$, but this is not universal and in many cases
one needs to perform an explicit calculation.  Luckily, for the
examples in this note, the representation-theoretic argument
will suffice.

\subsection{Equivariance of the cone construction}
\label{sec:equiv}

In order to calculate or simply identify the Killing superalgebras it
is often convenient to work in the cone.  This requires understanding
how to lift the calculation of the Lie derivative of a Killing
spinor along a Killing vector to the cone.  In \cite{JMFKilling} it is
shown that the cone construction is equivariant under the action of
the isometry group of $(M,g)$.  We will work at the level of the Lie
algebra.  Every Killing vector on $(M,g)$ defines a Killing vector on
the cone $(\overline M, \overline g)$.  Generically there are no other
Killing vectors on the cone, except in the case when $(M,g)$ is the
round sphere and hence the cone is flat. Let $X$ be a Killing vector
on $(M,g)$ and let $\overline X$ denote its lift to a Killing vector
on the cone.  Similarly let $\varepsilon$ be a Killing spinor on
$(M,g)$ and let $\overline\varepsilon$ denote the parallel spinor on
the cone to which it lifts.  Then it is proved in \cite{JMFKilling}
that
\begin{equation*}
  \eL_{\overline X} \overline\varepsilon = \overline{\eL_X
    \varepsilon}~,
\end{equation*}
which suggests a way to calculate the bracket $[-,-]:\fk_0 \otimes
\fk_1 \to \fk_1$:
\begin{itemize}
\item we lift the Killing vectors in $\fk_0$ and the Killing spinors
  in $\fk_1$ to Killing vectors and parallel spinors, respectively, on
  the cone;
\item we compute the spinorial Lie derivative there; and
\item we restrict the result to a Killing spinor on $(M,g)$.
\end{itemize}
Although somewhat circuitous, this procedure has the added benefit
that the Lie derivative of a parallel spinor is an algebraic
operation:
\begin{equation*}
  \eL_{\overline X} \overline\varepsilon = \rho(A_{\overline X})
  \overline\varepsilon~.
\end{equation*}
Since parallel spinors are determined by their value at any one point,
we can work at a point and we see that the above formula corresponds
to the restriction of the spin representation of $\fso_{n+1}$ to the
subalgebra corresponding to the image of $\fk_0$ in $\fso_{n+1}$,
acting on the subspace of the spinor module which is invariant under
the holonomy algebra of the cone.  For the case of the round spheres
which will occupy us in this paper, the holonomy algebra is trivial
and the isometries act linearly in the cone, whence $A_{\overline X} =
- \overline\nabla \overline X$ is actually constant.  Therefore the
above action is precisely the standard action of $\fk_0 = \fso_{n+1}$
on the relevant spinor module.

There is no need to lift the bracket $\fk_1  \otimes \fk_1 \to \fk_0$
to the cone, but it is possible to do this as well.  The only point to
notice is that in the cone we do not square parallel spinors to
parallel vectors, but to parallel 2-forms, which are constructed out
of the lifts of the Killing vectors on $(M,g)$.

\section{The Killing superalgebras of $S^7$, $S^8$ and $S^{15}$}
\label{sec:S7S8S15}

In this section we will exhibit the Killing superalgebras of some
low-dimensional spheres $S^n$, for $n=7,8,15$, and will show that they
are Lie algebras isomorphic to $\fso_9$, $\ff_4$ and $\fe_8$,
respectively.  The strategy is to exploit the equivariance of the cone
construction to show that these Killing algebras are isomorphic to the
Lie algebras constructed in \cite{AdamsExceptional}.

\subsection{$\fk(S^7) \cong \fso_9$}
\label{sec:S7}

The isometry Lie algebra of the unit sphere in $\RR^8$ is $\fso_8$,
acting via linear vector fields on $\RR^8$ which are tangent to the
sphere.  The $7$-sphere admits the maximal number of Killing spinors
of either sign of the Killing constant, which here is $8$.  Lifting
them to the cone, we have $\fso_8$ acting on the positive chirality
spinor module $\fS_+$ which is real and eight-dimensional.  The
Killing superalgebra is thus $\fk = \fso_8 \oplus \fS_+$ with the
following brackets: $\fso_8 \subset \fk$ is a Lie subalgebra, $\fso_8
\otimes \fS_+ \to \fS_+$ is the standard action and the map $\Lambda^2
\fS_+ \to \fso_8$ be the transpose of the previous map relative to the
inner products on both vectors and spinors.  The map is skewsymmetric
as shown because the spinor inner product is symmetric.  Therefore we
will obtain a Lie algebra.  Observe that triality says that
$\Lambda^2\fS_+ \cong \Lambda^2 V$, so that this map is actually an
isomorphism in this case.  The Jacobi identity requires the vanishing
of a trilinear map $\Lambda^3\fk \to \fk$.  The only component which
is in question is the one in $\Lambda^3\fS_+ \to \fS_+$.  Using the
inner product on $\fS_+$ we may identify this with an
$\fso_8$-invariant element in $\fS_+ \otimes \Lambda^3 \fS_+$, but it
may be shown the only such element is the zero map.  Indeed, letting
$\fS_+$, $\fS_-$ and $V$ have Dynkin indices $[0001]$, $[0010]$ and
$[1000]$, respectively, we find that $\Lambda^3\fS_+$ is irreducible
with Dynkin index $[1010]$, corresponding to the $56$-dimensional
kernel of the Clifford multiplication $V \otimes \fS_- \to \fS_+$.
Finally, a roots-and-weights calculation shows that
\begin{equation*}
  \fS_+ \otimes \Lambda^3\fS_+ \cong [0020] \oplus [0100] \oplus
  [1011] \oplus [2000]~,
\end{equation*}
whence there is no nontrivial invariant subspace.  The Lie algebra
structure just defined on $\fk$ is $36$-dimensional and coincides with
$\fso_9$.

\subsection{$\fk(S^8) \cong \ff_4$}
\label{sec:S8}

The isometry Lie algebra of the unit sphere in $\RR^9$ is $\fso_9$,
acting via linear vector fields on $\RR^9$ which are tangent to the
sphere.  The $8$-sphere admits the maximal number of Killing spinors
of either sign of the Killing constant, which here is $16$.  Lifting
them to the cone, we have $\fso_9$ acting on the spinor module $\fS$
which is real and sixteen-dimensional.  The Killing superalgebra is
$\fk = \fso_9 \oplus \fS$ with the following brackets: $\fso_9$ is a
Lie subalgebra, $\fso_9 \otimes \fS \to \fS$ is the standard action of
$\fso_9$ on its spinor representation, and $\Lambda^2\fS \to \fso_9$
to be the transpose of the standard action using the inner products on
vectors and spinors.  Since the spinor inner product is symmetric, the
map is skewsymmetric as shown.  This means that we will obtain a Lie
algebra.  The only nontrivial component of the Jacobi identity lives
in the subspace of $\fso_9$-equivariant maps $\Lambda^2\fS \to \fS$,
or using the inner product, an $\fso_9$-invariant element of $\fS
\otimes \Lambda^3\fS$.  However one can check that there are no such
invariants.  Indeed, since $\fS$ has Dynkin index $[0001]$, a
roots-and-weights calculation shows that
\begin{equation}
  \Lambda^3 \fS \cong [0101] \oplus [1001]~,
\end{equation}
where the representations on the right-hand side have dimensions $432$
and $128$, respectively.  Indeed, $[1001]$ is the kernel of the
Clifford multiplication $V \otimes \fS \to \fS$.  Tensoring the
first with $\fS$ we obtain
\begin{multline*}
  [0101] \otimes [0001] \cong [0002] \oplus [0010] \oplus [0100]
  \oplus [0102]\\
  \oplus [0110] \oplus [0200] \oplus [1002] \oplus [1010] \oplus
  [1100]~,
\end{multline*}
whereas tensoring the second with $\fS$ we obtain
\begin{equation*}
  [1001] \otimes [0001] = [0002] \oplus [0010] \oplus [0100] \oplus
  [1000] \oplus [1002] \oplus  [1010] \oplus [1100] \oplus [2000]~.
\end{equation*}
It is plain that there are no invariants in either expression.  The
resulting Lie algebra has dimension $36 + 16 = 52$ and can be shown
\cite{AdamsExceptional} to be a compact real form of $\ff_4$.  Unlike
the case of $\fso_9$ in Section~\ref{sec:S7}, here $\Lambda^2\fS \to
\fso_9$ is not an isomorphism: indeed $\Lambda^2\fS \cong \Lambda^2V
\oplus \Lambda^3V$.

\subsection{$\fk(S^{15}) \cong \fe_8$}
\label{sec:S15}

The isometry Lie algebra of the unit sphere in $\RR^{16}$ is
$\fso_{16}$.  The $15$-sphere admits the maximal number of Killing
spinors of either sign of the Killing constant, which here is $128$.
Lifting them to the cone, we have $\fso_{16}$ acting on the spinor
module $\fS_+$ which is real and $128$-dimensional.  The Killing
superalgebra is $\fk = \fso_{16} \oplus \fS_+$ with the following
brackets: $\Lambda^2\fso_{16} \to \fso_{16}$ is the Lie bracket,
$\fso_{16} \otimes \fS_+ \to \fS_+$ is the action of $\fso_{16}$ on
its half-spinor representation and $\Lambda^2\fS_+ \to \fso_{16}$ the
transpose map using the inner products.  As before, since the spinor
inner product is symmetric, the map is skewsymmetric as shown.  This
means that we will obtain a Lie algebra.  The resulting bracket can be
seen to satisfy the Jacobi identity.  Indeed, the only nontrivial
component of the Jacobi identity defines an $\fso_{16}$-equivariant
map $\Lambda^3\fS_+ \to \fS_+$.  Since the inner product is
non-degenerate on $\fS_+$, we can think of this as an
$\fso_{16}$-invariant element of $\fS_+ \otimes \Lambda^3\fS_+$, but
we can see that no such nontrivial element exists.  Indeed, letting
$[00000001]$ denote the Dynkin index of $\fS_+$, we find that
\begin{equation*}
  \Lambda^3\fS_+ \cong [00001001] \oplus [01000010] \oplus [10000001]~,
\end{equation*}
whence tensoring each of the modules in the right-hand side with $\fS_+$ we obtain
\begin{align*}
  [00001001] \otimes [00000001] = {} & [00000011] \oplus [00001000]
  \oplus [00001002] \oplus [00001100] \\
  & {} \oplus [00010011] \oplus [00011000] \oplus [00100002] \oplus [00100100] \\
  & {} \oplus [01000011] \oplus [01001000] \oplus [10000002] \oplus [10000100]~,
\end{align*}
\begin{align*}
  [00000001] \otimes [01000010] = {} & [00000011] \oplus [00001000]
  \oplus [00100000] \oplus [01000011]\\
  & {} \oplus [01001000] \oplus [01100000] \oplus [10000020] \oplus [10000100]\\
  & {} \oplus [10010000] \oplus  [11000000]~,
\end{align*}
and
\begin{align*}
  [00000001] \otimes [10000001] = {} & [00000011] \oplus [00001000]
  \oplus [00100000] \oplus [10000000]\\
  & {} \oplus [10000002] \oplus [10000100] \oplus  [10010000] \oplus [11000000]~.
\end{align*}
In all cases we see that there is nonzero invariant element.  The
resulting Lie algebra has dimension $120+128 = 248$ and can be shown
\cite{AdamsExceptional, GSW} to be isomorphic to the compact real form
of $\fe_8$.  Choosing $i\fS_+$ instead of $\fS_+$, we obtain the
maximally split real form of $\fe_8$ which has been the focus of
recent attention \cite{BBCE8}.  Notice that again $\Lambda^2\fS_+ \to
\fso_{16}$ is not an isomorphism, instead $\Lambda^2 \fS_+ \cong
\Lambda^2V \oplus \Lambda^6V$.

This construction of $\fe_8$ is also explained in \cite[§6.A]{GSW}, where
the nontrivial component of the Jacobi identity is proved
combinatorially using Fierz identities.

\section{Conclusion}
\label{sec:conc}

We have seen that a notion arising from supergravity, namely the
Killing superalgebra, when applied in a classical context, yields a
geometric construction of the exceptional Lie algebras of type $F_4$
and $E_8$.  This was accomplished by using Bär's cone construction to
relate the Killing superalgebra to the well-known construction of
these algebras using spin groups and their spinor representations.

There a number of things left to explore in relation to the
construction presented in this paper, some of which we are actively
considering:
\begin{itemize}
\item \emph{Further riemannian examples?}  The three examples
  considered here are of the following general form: $\fk = \fk_0
  \oplus \fk_1$ where $\fk_0$ is a Lie subalgebra, $\fk_1$ an
  $\fk_0$-module, there are $\fk_0$-invariant positive-definite inner
  products on $\fk_0$ and $\fk_1$ and hence on $\fk$ by declaring
  $\fk_0$ and $\fk_1$ to be orthogonal.  The bracket $\Lambda^2\fk_1
  \to \fk_0$ is defined precisely by the condition that the resulting
  inner product on $\fk$ be ad-invariant.  All but one components of
  the Jacobi identity of $\fk$ vanish. If Jacobi is satisfied, then we
  obtain a Lie algebra with a symmetric split and a positive-definite
  ad-invariant scalar product.  This means we have a riemannian
  symmetric space and in fact the nontrivial Jacobi identity is the
  algebraic Bianchi identity for the would-be curvature tensor.  We
  may therefore read off the possible such constructions from the list
  of symmetric spaces whose isotropy representation is spinorial, in
  which case the only examples are the above ones and the ones
  involving the exceptional Lie algebras $E_6$ and $E_7$, about which
  more below.  At any rate, we have looked explicitly at riemannian
  spheres in dimension $\leq 40$ which could give rise to Lie
  algebras, and have checked that the nontrivial Jacobi identity
  cannot follow trivially from representation theory.  It is therefore
  doubtful that other examples exist of \emph{precisely} this
  construction in riemannian signature.
\item \emph{Killing superalgebras of ``spheres'' in arbitrary
    signature.}  Considering other signatures (and hence
  possibly also imaginary Killing spinors) might provide geometric
  realisations of Lie superalgebras.
\item \emph{A similar construction for the remaining exceptional Lie
    algebras.}  In the case of $E_6$ and $E_7$, $\fk_0$ also contains
  ``R-symmetries'' which do not act geometrically on the manifold.
  Understanding these cases should help to understand conformal
  Killing superalgebras.  There does not seem to be a construction of
  $G_2$ using only spinors.
\item \emph{Of which structure on $S^{15}$ is $E_8$ the automorphism
    group?}  The existence of a Lie group is most naturally explained
  as automorphisms of some structure.  The construction of $E_8$ out
  of the $15$-sphere suggests that there ought to be some structure on
  $S^{15}$ of which $E_8$ is the automorphism group.  This may also
  provide a simple proof of the Jacobi identity without resorting to
  Fierz or roots-and-weights combinatorics.
\end{itemize}
I hope to report answers to some of these questions in the near
future.

\section*{Acknowledgments}

I have benefited from giving several talks on this topic at the
Dipartamento di Matematica ``U. Dini'' dell'Università degli Studi di
Firenze, at the Departamento de Análisis Matemático de la Universidad
de Alicante, and at the 18th North British Mathematical Physics
Seminar held at the University of York.  I am grateful to Dmitri
Alekseevsky and Andrea Spiro for arranging the visit to Florence, to
Salvador Segura Gomis for arranging the one to Alicante, and to Niall
MacKay for organising the meeting in York.  My interest in the Killing
superalgebra has been nurtured through collaboration with a number of
people, most recently, Emily Hackett-Jones, Patrick Meessen, George
Moutsopoulos, Simon Philip and Hannu Rajaniemi, to whom I offer my
thanks.  Finally, the roots-and-weights calculations in
Section~\ref{sec:S7S8S15} were performed using LiE \cite{LiE}, a
computer algebra package for Lie group computations.

\bibliographystyle{utphys}
\bibliography{AdS,AdS3,ESYM,Sugra,Geometry,Algebra}

\end{document}